\documentclass[12pt]{amsart}
\usepackage{amssymb}
\usepackage{amscd}

\newcommand{\Z}{{\mathbb Z}}
\newcommand{\R}{{\mathbb R}}
\newcommand{\F}{{\mathbb F}}
\newcommand{\Line}{{\mathbb P}^1}

\newcommand{\Fq}{{\F_q}}

\newcommand{\bP}{{\mathbb P}}
\newcommand{\A}{{\mathbb A}}

\newcommand{\gt}{{\widetilde{g}}}
\newcommand{\Nqab}{{N^{\rm ab}_q(g)}}

\newcommand{\Nqtor}{{N^{\rm tor}_q(g)}}

\newcommand{\Nqtc}{{N^{\rm tc}_q(g)}}

\DeclareMathOperator{\Hom}{Hom}
\DeclareMathOperator{\Spec}{Spec}

\theoremstyle{plain}
\newtheorem  {thm}        {Theorem}     [section]
\newtheorem  {lemma}[thm] {Lemma}
\newtheorem  {prop} [thm] {Proposition}

\newtheorem* {thm*}       {Theorem}
\newtheorem* {prop*}      {Proposition}

\theoremstyle{definition}

\theoremstyle{remark}
\newtheorem* {remark}     {Remark}

\begin{document}


\title[Curves of every genus]{
  Curves of every genus with many points, I: Abelian and toric families
}

\author[A. Kresch]{Andrew Kresch}
\address{
  Department of Mathematics,
  University of Pennsylvania,
  Philadelphia, PA 19104--6395
}
\email{kresch@math.upenn.edu}

\author[J. Wetherell]{Joseph~L.~Wetherell}
\address{
  Mathematical Sciences Research Institute,
  Berkeley, CA 94720-5070
}
\email{jlwether@alum.mit.edu}

\author[M. Zieve]{Michael~E.~Zieve}
\address{
  Center for Communications Research,
  Princeton, NJ 08540--3699
}
\email{zieve@idaccr.org}


\thanks{The authors thank Andrew Granville, David Harbater,
Sandy Kutin, Hendrik W. Lenstra, Jr., Jean-Pierre Serre, Chris Skinner, 
William Stein, Jeff VanderKam, and Hui Zhu for
helpful comments and conversations.  The authors were supported
in part by NSF Mathematical Sciences Postdoctoral Research Fellowships.
}

\begin{abstract}
Let $N_q(g)$ denote the maximal number of $\Fq$-rational
points on any curve of genus $g$ over $\Fq$.  Ihara (for square $q$)
and Serre (for general $q$) proved that
$\limsup_{g\to\infty} N_q(g)/g > 0$ for any fixed $q$.
Here we prove
$\lim_{g\to\infty} N_q(g)=\infty$.  More precisely, we use
abelian covers of $\Line$
to prove $\liminf_{g\to\infty} N_q(g)/(g/\log g)>0$, and we use
curves on toric surfaces to prove $\liminf_{g\to\infty} N_q(g)/g^{1/3}>0$;
we also show that these results are the best possible that can be proved
using these families of curves.
\end{abstract}


\maketitle


\section{Introduction}

Let $\Fq$ be the field with $q$ elements, and let $C$ be a curve
(nonsingular, projective, geometrically irreducible) of
genus $g$ defined over $\Fq$.
The Riemann hypothesis for curves over finite fields (Weil's theorem)
implies that the number of $\Fq$-rational points on $C$ satisfies
the inequality
$\#C(\Fq)\leq q+1+2g\sqrt{q}$.  Ihara~\cite{Ih} showed that, when
$g$ is large compared to $q$, this inequality can be improved
significantly.  In this paper we study such improvements, seeking
the best upper bound on $\#C(\Fq)$ which depends only on $g$ and
$q$ -- we let $N_q(g)$ denote this bound, i.e., $N_q(g)$ is the
maximum number of $\Fq$-rational points on any curve over $\Fq$
of genus $g$.

The Weil bound implies that, for $q$ fixed and $g$ varying,
$N_q(g)\leq g(2\sqrt{q}) + O_q(1)$.
Ihara~\cite{Ih} observed that equality cannot hold when $g$ is much
larger than $q$, since this would
imply the existence of curves having a negative number of points over
$\F_{q^2}$.  This observation was extended by Drinfeld and
Vladut~\cite{DV} to derive the bound
$N_q(g)\leq g(\sqrt{q}-1)+o_q(g)$, or in other words (for $q$ fixed)
$\limsup_{g\to\infty} N_q(g)/g \leq \sqrt{q}-1$.

In order to test the sharpness of the Drinfeld-Vladut bound,
it is necessary to produce curves with many points.
Five essentially different approaches have been used.
Serre~\cite{Se,Se4} used class field towers to show that,
for any $q$, we have
$\limsup_{g\to\infty} N_q(g)/g\geq \gamma_q>0$; subsequently, other
authors have used the same method to derive the same result but with
larger constants $\gamma_q$ (see~\cite{NX} and the references therein).
Ihara~\cite{Ih-66,Ih-69,Ih-75,Ih-79,Ih} used supersingular
points on Shimura curves to show that, when $q$ is a square, one can
take $\gamma_q=\sqrt{q}-1$ -- the largest constant possible, according
to the Drinfeld-Vladut result; subsequently, Manin and Vladut~\cite{MV}
used supersingular points on Drinfeld modular curves to derive the same result.
(Some special cases of Ihara's construction were rediscovered in~\cite{TVZ}.)
Garcia and Stichtenoth wrote down explicit towers of Artin-Schreier
extensions~\cite{GS1,GS3} and (jointly with Thomas) Kummer
extensions~\cite{GST} which have many points (interestingly, Elkies has shown
that several of the Garcia-Stichtenoth towers are examples of towers of
modular~\cite{El98} or Drinfeld modular~\cite{El3ecm} curves).
Zink~\cite{Zi} showed that certain
degenerate Shimura varieties are curves with many
points of degree three over the prime field.
The fifth and most recent approach is that of Frey, Kani and
V\"olklein~\cite{FKV}, who combine rigidity methods from group theory
with a careful analysis of certain abelian varieties in order to
produce curves over $\Fq$ having unramified covers of
arbitrarily large degree in which some $\Fq$-point splits completely.

The above results exhibit sequences of genera $g$ for which
$N_q(g)$ is `large' (relative to $g$).  In the present paper we
examine how small $N_q(g)$ can be (relative to $g$).
Our first new result is

\begin{thm}
\label{thm-infinite}
For fixed $q$, we have $\lim_{g\to\infty} N_q(g) = \infty$.
\end{thm}

The difficulty in proving this result is that we insist on finding
curves in {\em every} large genus.  The methods listed above for producing
curves with many points do not enable one to do this (for example, modular
curves achieve few genera~\cite{CWZ}).
On the other hand, the simplest families of curves which do attain
every large genus -- hyperelliptic, trigonal, bielliptic --
are all low-degree covers of low-genus curves, so they
cannot have many points.  For instance, no curve of these three types can have
more than 10 points over $\F_2$, and in fact before the present paper it
was not known whether $N_2(g)>10$ for all large $g$.

We will give three proofs of Theorem~\ref{thm-infinite}, based on
studying three particular family of curves:
tame cyclic covers of $\Line$,
arbitrary abelian covers of $\Line$,
and curves embedded in toric surfaces.
In each case we prove a lower bound on $N_q(g)$, and in the
abelian and toric cases, we
show that (up to a universal constant factor) these are the best lower
bounds on $N_q(g)$ provable with these families of curves.
For each $q$ and $g$, let $\Nqab$ denote the maximum number of $\Fq$-rational
points on any genus-$g$ curve which is an abelian cover of $\Line$ over $\Fq$;
let $\Nqtc$ and $\Nqtor$ denote the corresponding quantities for tame cyclic
covers of $\Line$ over $\Fq$ and for curves which embed in toric surfaces
over $\Fq$, respectively.

\begin{thm}
\label{thm-main-abelian}
For any fixed $q$, there exist constants $0<a_q<b_q$ such that:
for every $g>1$, we have
$ a_q\cdot g/\log{g} < \Nqab < b_q\cdot g/\log{g}.$
\end{thm}

\begin{thm}
\label{thm-main-toric}
For any fixed $q$, there exist constants $0<c_q<d_q$ such that:
for every $g>0$, we have
$ c_q\cdot g^{1/3} < \Nqtor < d_q\cdot g^{1/3}.$
\end{thm}

\begin{thm}
\label{thm-main-cyclic}
For any fixed $q$, there exist constants $e_q,f_q>0$ such that:
for every $g>f_q$, we have
$ e_q \cdot g(\log\log g)/(\log g)^3 < \Nqtc.$
\end{thm}


These three theorems are proved in Sections \ref{sec-abelian}, \ref{sec-toric},
and \ref{sec-tame}, respectively.
The upper bound in Theorem~\ref{thm-main-abelian} is due to
Frey, Perret, and Stich\-te\-noth~\cite{FPS}.  The other
inequalities are new.

Several authors have previously used abelian covers
to produce curves with many points.
Serre~\cite{Se,Se4} proposed abelian covers
of known curves as a convenient source of curves with many points in
case $g$ is not much larger than $q$, the idea being that such covers can be
understood via class field theory.
Other authors have subsequently taken this approach to produce numerous
examples (cf.~\cite{GV} or \cite{NX0}).
The hard part in our work is finding examples with prescribed genus.
The toric approach of Section~\ref{sec-toric} is new.

Theorems \ref{thm-main-abelian} and \ref{thm-main-toric} imply that any one
genus behaves in roughly the same manner as any other, with respect to the
maximum number of $\Fq$-rational points on curves of this genus lying in
either of two special families of curves.  In subsequent work, we have shown
a similar assertion for the family of {\em all} curves:

\begin{thm}[\cite{phewkz}]
For any fixed $q$, there exists $h_q>0$ such that:
for every $g>0$, we have
$h_q\cdot g < N_q(g).$
\end{thm}

Throughout this paper, by a curve over a field $k$, or simply
`curve', we mean a complete non-singular one-dimensional variety
defined over $k$ which is geometrically irreducible.  We reserve the symbol
$p$ for the characteristic of the field under consideration.

We advise the reader that the sections of this paper are logically
independent from each other, and can be read in any order.


\section{Abelian covers of $\Line$}
\label{sec-abelian}

In this section we prove Theorem~\ref{thm-main-abelian} in the following form,
where $\Nqab$ denotes the maximum number of $\Fq$-rational points on any
genus-$g$ curve which admits an abelian cover of $\Line$ over $\Fq$:

\begin{thm}
\label{thm-wild}
For any fixed $q$, we have
$$
\inf_{g>1}\frac{\Nqab}{g/\log g} > 0
\quad\quad \text{and} \quad\quad
\sup_{g>1}\frac{\Nqab}{g/\log g}<\infty.
$$
\end{thm}

The upper bound is due to Frey, Perret and Stichtenoth~\cite{FPS}.
In this section we prove the lower bound.  The specific abelian covers
we use will be fiber products of elementary abelian covers with a single
degree-2 cover.

First consider elementary abelian $p$-covers of $\Line/\Fq$ which are only
ramified at $x=0$.  (Here $p$ denotes the characteristic of $\Fq$.)
For instance, consider the equations
\begin{equation*}
    y_0^p - y_0 = x^{-i_0}, \quad
    y_1^p - y_1 = x^{-i_1}, \quad
                \dots, \quad
    y_{n-1}^p - y_{n-1} = x^{-i_{n-1}},
\end{equation*}
where $i_0<i_1<\dots<i_{n-1}$ is an increasing sequence of positive
integers coprime to $p$.  It can be shown that these equations define a
curve of genus
$$
 \frac{p-1}{2}\Bigl((i_0-1) + (i_1-1)p + \dots + (i_{n-1}-1)p^{n-1} \Bigr).
$$
For $p=2$ and any fixed $n$, these curves achieve every sufficiently large
genus; this is the crux of our proof of Theorem~\ref{thm-wild}
for even $q$ (which we give at the end of this section).
For odd $p$, however, the genus of such a curve is divisible by $(p-1)/2$ --
a problem we will solve by taking a degree-2 cover of our curve.  A more
serious problem is that the genus is never congruent to $(p+1)/2$ mod~$p$;
to attain every congruence class, we allow ramification at two points.

\begin{lemma}
\label{wild-equations}
Let $i_0<i_1<\dots<i_{n-1}$ and $j_0<j_1<\dots<j_{n-1}$ be
increasing sequences of positive integers which are coprime to $p$.
Then the equations
\begin{align*}
    y_0^p - y_0 &= x^{-i_0}(x-1)^{-j_0} \\
    y_1^p - y_1 &= x^{-i_1}(x-1)^{-j_1} \\
                &\;\;\vdots \\
    y_{n-1}^p - y_{n-1} &= x^{-i_{n-1}}(x-1)^{-j_{n-1}}
\end{align*}
describe a curve $C/\Fq$ of genus
$$\frac{p-1}{2}
\Bigl((i_0+j_0) + (i_1+j_1)p + \dots + (i_{n-1}+j_{n-1})p^{n-1} \Bigr).$$
The map $C\to\Line_x$ is Galois with Galois group $(\Z/p\Z)^{n}$.
It is unramified away from $x=0$ and $x=1$; furthermore, $x=\infty$
splits completely in the cover $C\to\Line$.
\end{lemma}

\begin{proof}
We use standard facts about composita of Artin-Schreier extensions, cf.\ \cite{GVJCT}
or \cite{GS0}.
Since the right-hand sides of the defining equations are linearly independent
over $\F_p$, the field $L:=\F_q(x, y_0, \ldots, y_{n-1})$ is a Galois extension
of $\F_q(x)$ with Galois group $(\Z/p\Z)^n$, and the genus of $L$ is
the sum of the genera of all the degree-$p$ subextensions of $L$.
The formula for the genus now follows from the easy fact that if $i$ and $j$
are positive integers coprime to $p$, and $f(x)\in \F_q[x]$ has
degree less than $i+j$ and is coprime to $x(x-1)$, then the
curve over $\F_q$ defined by
$$y^p-y=f(x)x^{-i}(x-1)^{-j}$$
has genus $((p-1)/2)(i+j)$.
The remaining assertions are clear.
\end{proof}

In our application of this lemma we will want sequences $i_\nu$ and $j_\nu$
giving a genus on the order of $np^n$, and moreover we require these
sequences to yield such genera in every residue class mod $p^n$.
We address these issues in the following combinatorial lemma.

\begin{lemma}
\label{congruence-problem}
Let $n$ and $d$ be positive integers.
If $p$ is an odd prime, then there exist
increasing sequences $i_0<\dots<i_{n-1}$ and $j_0<\dots<j_{n-1}$ of
positive integers coprime to $p$ such that each $i_k+j_k<(p+3)(k+1)$
(for $0\leq k\leq n-1$) and
$\sum_{k=0}^{n-1} (i_k+j_k)p^k\equiv d$ {\rm (mod}~$p^n)$.
\end{lemma}

\begin{proof}
We prove the lemma by induction on $n$.  Throughout the proof,
if $r$ is an integer, then
$r_p$ denotes the unique integer such that $0 \le r_p \le p-1$ and
$r_p \equiv r$ (mod~$p$).  For $n=1$ we must find positive integers
$i_0$ and $j_0$ which are coprime to $p$ such that
$d\equiv i_0+j_0$ (mod~$p$) and $i_0+j_0 < p+3$.
This is easily done: if $d_p>1$, then set
$i_0=d_p-1$ and $j_0=1$; and if
$d_p\le1$, then set $i_0=d_p+1$ and $j_0=p-1$.

Now assume, inductively, that we have sequences $i_0<\dots<i_{n-2}$
and $j_0<\dots<j_{n-2}$ of positive integers coprime to $p$ such that
each $i_k+j_k<(p+3)(k+1)$ and $S:=\sum_{k=0}^{n-2}(i_k+j_k)p^k$ is congruent
to $d$~(mod~$p^{n-1}$).  We will find integers $i_{n-1}$ and $j_{n-1}$ which
are
coprime to $p$, where $i_{n-2}<i_{n-1}$ and $j_{n-2}<j_{n-1}$, such that
$i_{n-1}+j_{n-1}<(p+3)n$ and $S+(i_{n-1}+j_{n-1})p^{n-1}\equiv d$ (mod~$p^n$).
This last congruence is equivalent to $i_{n-1}+j_{n-1}\equiv d'$ (mod~$p$),
where $d'=(d-S)p^{1-n}$.
Now we let
$$i_{n-1}=i_{n-2}+b  \quad\text{and}\quad
  j_{n-1}=j_{n-2}+1+(d'-1-b-i_{n-2}-j_{n-2})_p,$$
where $b\in\{1,2,3\}$ is chosen so that $i_{n-1}$ and $j_{n-1}$ are
both coprime to $p$.
Then the desired conditions on $i_{n-1}$ and $j_{n-1}$ are satisfied,
and we have completed the induction.
\end{proof}

We now prove Theorem~\ref{thm-wild}, first for odd $q$ and then for even $q$.

\begin{proof}[Proof of Theorem~\ref{thm-wild} for odd $q$]
Fix a finite field $\Fq$ of characteristic $p>2$.  Let $f(x)\in\Fq[x]$ be a monic
irreducible polynomial of even degree.
For any $n>0$, consider the system of equations
\begin{align*}
    y_0^p - y_0 &= x^{-i_0}(x-1)^{-j_0} \\
    y_1^p - y_1 &= x^{-i_1}(x-1)^{-j_1} \\
                &\;\;\vdots \\
    y_{n-1}^p - y_{n-1} &= x^{-i_{n-1}}(x-1)^{-j_{n-1}} \\
                    y^2 &= f(x),
\end{align*}
where both $i_0<i_1<\dots<i_{n-1}$ and $j_0<j_1<\dots<j_{n-1}$ are
increasing sequences of positive integers coprime to $p$.
This is the fiber product of a degree-$p^n$ cover
$\phi\colon C\to\Line_x$ (as discussed in Lemma~\ref{wild-equations})
with a degree-$2$ cover $\psi$ from the curve $y^2=f(x)$ to $\Line_x$.
Since the two covers have coprime degrees, the system of
equations describes a curve $\widetilde{C}$ over $\Fq$.  Moreover,
the induced degree-$2p^n$ cover $\widetilde{C}\to\Line_x$ is abelian,
since it is the fiber product of abelian covers.  Write the
degree of the different of $\psi$ as $2D$.
Then by the Riemann-Hurwitz formula applied to $\widetilde{C}\to C$,
the genus $\widetilde{g}$ of $\widetilde{C}$ is given by
\begin{equation*}
\widetilde{g} = (p-1)\bigl((i_0+j_0)+(i_1+j_1)p+\dots+(i_{n-1}+j_{n-1})p^{n-1}\bigr) + p^n D - 1.
\end{equation*}
Since $f$ is monic and has even degree,
$x=\infty$ splits completely under the map $\widetilde{C}\to\Line_x$,
so $\widetilde{C}$ has at least $2p^n$ rational points over $\Fq$.

Pick any $g>p^2+3p$; we now describe choices of the parameters above so
that $\widetilde{g}=g$.  
Let $n$ be the largest integer such that $(p+3)np^n<g$.  Our assumption
implies $n \ge 1$.  Lemma~\ref{congruence-problem} yields sequences
$i_0<i_1<\dots<i_{n-1}$ and $j_0<j_1<\dots<j_{n-1}$ of positive integers
coprime to $p$ such that each $i_k+j_k<(p+3)(k+1)$ and
$\sum_{k=0}^{n-1}(i_k+j_k)p^k$ is congruent to $(g+1)/(p-1)$ mod~$p^n$.
Define $D$ by the equation
\begin{equation*}
g = (p-1)\sum (i_k+j_k)p^k + p^nD - 1.
\end{equation*}
Then $D$ is an integer, and our bound on $i_k+j_k$ implies $D$ is positive.
%
%
Let $f(x)$ be a monic irreducible polynomial in $\Fq[x]$
of degree $2D$; then the degree of the different of the degree-$2$ cover
from $y^2=f(x)$ to $\Line_x$ is $2D$.  Now $f$, $i_k$, and $j_k$ satisfy
all the conditions of the previous paragraph, and yield a curve
$\widetilde{C}$ of genus $g$ such that $\#\widetilde{C}(\Fq)\geq 2p^n$.
Our choice of $n$ implies that $\log g > n\log p$ and
$p^n\geq g/(p(p+3)(n+1))$, so
\begin{equation*}
   \#\widetilde{C}(\Fq) \cdot \log g >
     2p^n\cdot n\log p \ge \frac{2 \log p}{(p+3)p}
     \frac{n}{(n+1)} g \ge \frac{\log p}{2p^2} g.
\end{equation*}

We have shown that $\inf_{g>p^2+3p}\Nqab/(g/\log g)$ is positive.  To complete
the proof, we must show that $\Nqab>0$ for all $g>1$.  For this, let $h(x)\in\Fq[x]$
be squarefree of degree $2g+1$, and note that
the curve $y^2=h(x)$ has genus $g$ and has an $\Fq$-rational point.
\end{proof}

%
%
%

\begin{proof}[Proof of Theorem~\ref{thm-wild} for even $q$]
Let $i_0<i_1<\dots<i_{n-1}$ be an increasing sequence of odd positive
integers, and consider the system of equations
\begin{equation*}
    y_0^2 + y_0 = x^{-i_0}, \quad
    y_1^2 + y_1 = x^{-i_1}, \quad
                \dots, \quad
    y_{n-1}^2 + y_{n-1} = x^{-i_{n-1}}.
\end{equation*}
Just as in the proof of Lemma~\ref{wild-equations}, we see that these
equations define a curve $C$ over $\F_2$ such that $\#C(\F_2)\geq 2^n$
and the genus of $C$ is
$$
\frac{i_0-1}{2}2^0+\frac{i_1-1}{2}2^1+\dots+\frac{i_{n-1}-1}{2}2^{n-1}.
$$
For a fixed positive integer $n$, we can find sequences $i_k$ as above
yielding curves of any genus greater than $n2^{n+1}-5$: there is a
unique choice of $i_0\in\{1,3\}$, $i_1\in\{5,7\}$, \dots,
$i_{n-2}\in\{4n-7,4n-5\}$ for which $\sum_{k=0}^{n-2}(i_k-1)2^{k-1}$
attains any prescribed integer value in the interval
$[(n-3)2^n+4,(n-3)2^n+2^{n-1}+3]$, and for any $g\geq n2^{n+1}-4$
there is then an odd integer $i_{n-1}>4n-5$ yielding a genus-$g$ curve.

For a given nonnegative integer $g$, let $n$ be the unique positive
integer such that $n2^{n+1}-4\leq g< (n+1)2^{n+2}-4$.
Then the previous paragraph shows that $\Nqab\geq 2^n$,
and for $g>1$ this implies that $\Nqab>((\log 2)/4)\cdot g/\log g$.
\end{proof}


\section{Curves in toric surfaces}
\label{sec-toric}

Curves in toric surfaces are those given by a single equation $f(x,y)=0$,
such that the curve in $\bP^2$ defined by $f$ has a resolution of
singularities of a specific form.
As we explain below, the genus and number of rational points of such a curve
are governed by the shape of the definining equation $f(x,y)$---specifically,
by the Newton polygon of $f$.
Letting $\Nqtor$ denote the maximum number of $\Fq$-rational points on
any genus-$g$ curve which embeds in a toric surface over $\Fq$, we will prove
Theorem~\ref{thm-main-toric} in the following form:

\begin{thm}
\label{thm-toric}
For any fixed $q$, we have
\begin{equation*}
\inf_{g>0}\frac{\Nqtor}{g^{1/3}} > 0
\quad\quad \text{and} \quad\quad
\sup_{g>0}\frac{\Nqtor}{g^{1/3}}<\infty.
\end{equation*}
\end{thm}

The following result presents the curves we use to prove the lower bound.
\begin{prop}
\label{genuscomputation}
Let $k$ be a field of characteristic $p>0$, and let $r$ be a positive integer.
Choose integers $0<a_0<\dots<a_r$.
Then the equation $f(x,y)=0$, where
\begin{equation*}
f(x,y):=1+y+x^{r+1}+\sum_{i=0}^r x^iy^{p(a_i+\cdots+a_r)},
\end{equation*}
defines a curve over $k$ with at least $r$ rational points and genus
\begin{equation}
\label{newtonpoints}
g:=-r+p\sum_{i=0}^r ia_i.
\end{equation}
Moreover, this curve has a smooth complete model in some smooth
projective toric surface over $k$.
\end{prop}

\begin{remark}
The {\em Newton polygon} of a polynomial $f\in k[x,y]$ is the
convex hull in $\R^2$ of the set of lattice points $(i,j)$ for
which the coefficient of $x^iy^j$ in $f$ is nonzero.
The expression for the genus (\ref{newtonpoints}) is the number
of interior lattice points in the Newton polygon of $f$.
\end{remark}

Proposition \ref{genuscomputation} is an immediate consequence
of the following algebro-geometric statement.
\begin{prop}
\label{agprop}
Let $k$ be a field.
Let $f\in k[x,y]$ satisfy
\begin{itemize}
\item[(i)] $f$, $\partial f/\partial x$, and $\partial f/\partial y$
generate the unit ideal in $k[x,y]$;
\item[(ii)] $f$ has nonzero constant term, and $f$ is not in $k[x]$
or $k[y]$;
\item[(iii)] every lattice point on the boundary of the Newton
polygon of $f$ is either a vertex or lies on the horizontal or vertical
coordinate axis.
\end{itemize}
Let $g$ be the number of lattice points in the interior of the
Newton polygon of $f$, and let $v$ be the number of vertices of
the Newton polygon.
Then $f(x,y)=0$ defines a curve which has genus $g$ and at least $v-2$
rational points over $k$.
Moreover, the curve admits a smooth complete model in some smooth
projective toric surface over $k$.
\end{prop}



The plan for the rest of this section is as follows.
Given Proposition \ref{agprop}, the lower bound of Theorem \ref{thm-toric}
follows by easy combinatorics, which we do first.
The upper bound uses known combinatorics of lattice polygons,
coupled with some geometry of curves in algebraic surfaces.
This bit of algebraic geometry is also what is needed to prove
Proposition \ref{agprop}, and we present this second.
Lastly, we establish the upper bound in Theorem \ref{thm-toric}.

\begin{proof}[Proof that Proposition~\ref{agprop} implies the lower bound of
Theorem~\ref{thm-toric}]
For any $g>0$, there is a genus-$g$ curve which has at least one
$\Fq$-rational point and which satisfies (i)--(iii) of Proposition~\ref{agprop}:
for instance, let $h(x)\in\Fq[x]$ be squarefree of degree $2g+1$, and take
$y^2+y=h(x)$ for even $q$ and $y^2=h(x)$ for odd $q$.
Thus, it suffices to show that for sufficiently large $g$ there exist
curves in toric surfaces whose genus is bounded above by a constant times
the cube of the number of rational points.
Moreover, it suffices to consider the case where $q$ is a prime $p$.

We use Proposition \ref{genuscomputation} (which follows at once from
Proposition~\ref{agprop}).
We may consider each residue class of $g$ mod~$p$ separately.
By dint of (\ref{newtonpoints}),
we are reduced to showing that for $r$ in a given residue class mod~$p$,
there are
increasing nonnegative sequences $\{a_i\}$ with $g/r^3$ bounded,
so that
the sums $\sum_{i=0}^r ia_i$ take on all sufficiently large positive integers.
Consider sequences $\{a_i\}$ with $a_0=0,\dots,a_{r-2}=r-2$ and
$r-1\le a_{r-1}\le 2(r-1)<a_r$.  Since nonnegative linear combinations of
$r-1$ and $r$ achieve all integer values greater than $r^2$, there are
sequences of this type (for fixed $r$) for which
$\sum_{i=0}^r ia_i$ takes on any prescribed value greater than $(2r^3+15r^2)/6$.
This proves the lower bound in Theorem \ref{thm-toric}.
\end{proof}

Now we prove Proposition \ref{agprop}.
Any $f(x,y)$ satisfying (i)--(iii) defines a variety on a suitable
toric surface having the claimed (arithmetic) genus by a result in \cite{h},
but it is not immediately clear that this variety is a curve
(i.e., smooth and geometrically irreducible).
So we need to show more: examination of its
defining equations in coordinate charts establishes smoothness,
and then a little intersection theory eliminates the possibility of
geometric reducibility.
Given the analysis we need for this, the genus computation falls out
easily; for more general results see \cite{h}.
The reader is only assumed to know basic facts about toric varieties,
particularly about toric surfaces; cf.~\cite{f}.

\begin{lemma}
\label{arithgenusformula}
Let $k$ be a field.
Let $f(x,y)\in k[x,y]$ be a nonzero polynomial, and let $\Gamma$ be the
Newton polygon of $f$.
Then the variety defined by
$f(x,y)=0$ can be compactified to a variety $V$ in a smooth
projective toric surface $Y$, such that $V$ is disjoint from the
set of fixed points for the toric action on $Y$.
Moreover, for any smooth projective toric surface $Y$ admitting such
a compactification $V$, the arithmetic genus of $V$ is equal to the
number of lattice points in the interior of $\Gamma$.
\end{lemma}

\begin{proof}
A toric surface is a two-dimensional normal variety $Y$,
equipped with the action of the two-dimensional algebraic torus
$T^2=(\A^1\smallsetminus\{0\})^2$ and a dense equivariant
embedding $T^2\to Y$.
The combinatorial object attached to $Y$ is a {\em fan}, which we
denote $\Delta$, consisting of cones in a two-dimensional lattice $N$;
the dual lattice $\Hom(N,\Z)$ (usually denoted $M$) is identified with
the lattice $\Z^2$ in which the Newton polygon sits (so $x$ and $y$ can
be regarded as coordinates on $Y$).
Cones of the fan correspond to affine coordinate charts on $Y$.

Each ray (one-dimensional cone) of $\Delta$ defines a family of
half-planes in $\Z^2$, namely $\rho=\R_{+}\cdot v$ defines
a half-plane $\{\,\alpha\,|\,\alpha(v)\ge h\,\}$ for any $h\in \R$.
Let us say that any half-plane in this family is
{\em associated with $\rho$}.
Consider the variety in $T^2$ defined by the equation $f(x,y)=0$,
with closure $V$ in $Y$.
We show now, by transforming $f$ into local coordinate systems,
that a complete toric surface $Y$
has all its fixed points disjoint from $V$
if and only if every supporting half-plane which meets $\Gamma$ in an edge
is associated with some ray $\rho\in\Delta$.
Indeed, consider a two-dimensional cone
$\nu=\R_{+}\cdot v + \R_{+}\cdot w$;
the corresponding affine chart $U_\nu$ has coordinate ring
$k[G_\nu]$, where $G_\nu$ is the semigroup
$$\{\,\alpha\in\Hom(N,\Z)
\,|\,\alpha(v)\ge 0 {\rm\ and\ }\alpha(w)\ge 0\,\}.$$
Now $V$ is disjoint from the origin of $U_\nu$ if and only if
there is an element in $k[G_\nu]$, with nonzero constant term,
which is equal to $f(x,y)x^ry^s$ for some integers $r$ and $s$.
Such an element exists if and only if some translate of $\Gamma$
is contained in $G_\nu$ and contains the origin, i.e.,
if the corner point of the intersection of
the pair of supporting half-planes of $\Gamma$ associated with $v$ and with $w$
is a vertex of $\Gamma$.

The nonsingular toric surfaces are those with $U_\nu\simeq \A^2$
for every two-dimensional cone $\nu\in\Delta$, or equivalently,
with every two-dimensional cone generated by
vectors which form a $\Z$-basis for $N$.
In this case we say $\Delta$ is {\em nonsingular}.
A toric surface is projective if and only if it is a complete variety,
and this is the case if and only
if the union of the cones in its fan is equal to $N$:
such a fan is called {\em complete}.

Any finite set of rays in $N$ is contained in a nonsingular complete
fan $\Delta$.
Consequently, given any nonzero polynomial $f(x,y)$, there exists a
nonsingular projective toric surface $Y$ such that
the subvariety $V$ of $Y$ defined by $f(x,y)=0$ (as above) is
disjoint from the fixed points of $Y$.
Now we consider the long exact sequence of
sheaf cohomology groups
\begin{equation}
\label{sheafcoho}
H^1(Y,{\mathcal O}_Y)\to
H^1(V,{\mathcal O}_V)\to
H^2(Y,{\mathcal O}_Y(-V))\to
H^2(Y,{\mathcal O}_Y).
\end{equation}
The first and last terms in (\ref{sheafcoho}) vanish because
they are invariant under blowing up a
rational point on a projective surface, and they vanish for $\bP^2$.
So the arithmetic genus of $V$ is equal to
$\dim H^2(Y,{\mathcal O}_Y(-V))$.
By Serre duality this equals
$\dim H^0(Y,K_Y(V))$, where
$K_Y$ is the canonical bundle of $Y$.
But $-K_Y$ is the sum of the toric divisors
(closures of one-dimensional torus orbits) of $Y$, so one can identify
the set of lattice points in the interior of $\Gamma$ with a basis
for $H^0(Y,K_Y(V))$.
\end{proof}

We continue with the notation of the Lemma -- toric surface $Y$
with subvariety $V$ disjoint from the fixed point set of $Y$ -- and we
describe the intersection of
$V$ with any of the one-dimensional torus orbits of $Y$.
Let $\rho=\R_+\cdot v$ be a ray of $\Delta$,
and consider the corresponding torus orbit $E$.
Corresponding to $\rho$ is a toric chart $U_\rho$, and we can identify
$U_\rho\simeq \Spec k[t,u,u^{-1}]$
so that $t=0$ defines $E$.
Transforming $f$ into $(t,u)$-coordinates and setting $t=0$
yields
\begin{equation}
\label{modpu}
V\cap E \simeq \Spec k[u, u^{-1}]/(p(u)),
\end{equation}
where $p(u)$ is the Laurent polynomial whose
sequence of coefficients, indexed by $\Z$, is equal to
the sequence of coefficients of monomials $x^iy^j$ of $f$,
for $(i,j)$ lying on the boundary $\ell$
of the half-plane supporting $\Gamma$,
associated with $\rho$
(the identification of $\Z$ with the set of lattice points on $\ell$
is to be via an affine linear map,
so $p(u)$ is defined only up to
multiplication by a power of $u$ and interchanging $u$ and $u^{-1}$).

If we understand the degree of the Laurent polynomial $p(u)$ to be
the maximal degree minus the minimal degree of all the monomials appearing
in $p(u)$, then
the degree of $p(u)$ is one less than the number
of lattice points in $\ell\cap \Gamma$.
In particular, if the $\ell$ associated with $\rho$
contains no lattice points in $\Gamma$ other than vertices of $\Gamma$,
then $p(u)$ has degree $0$ or $1$.
Hence the intersection $V\cap E$ is either empty or consists of a single
$k$-valued point, which is a regular point of $V$.

Now suppose $f(x,y)$ satisfies (i)--(iii) of Proposition \ref{agprop}.
Since, by (ii), $\Gamma$ is contained in the first quadrant of $\Z^2$ and
contains the origin, there exists a smooth projective toric surface
$Y$ with fixed point set disjoint from $V$, such that $Y$ is obtained
by starting with $\bP^2$ and repeatedly blowing up points in
the complement of $\A^2\subset \bP^2$.
Now, then, the test for $V\cap \A^2$ to be nonsingular
is precisely condition (i).
Additionally,
every point of $V$ which is not in $\A^2$ is a regular point of $V$ by
condition (iii) and the previous paragraph.
Hence $V$ is nonsingular.
The genus assertion is Lemma \ref{arithgenusformula}, and $V$ has at
least $v-2$ rational points, one from each of $v-2$
intersections (\ref{modpu}) with $\deg p(u)=1$
(where $v$ is the number of vertices of the Newton polygon).
It remains only to show that $V$ is absolutely irreducible.

We pass to the algebraic closure of $k$, and we suppose
$f=f_1f_2$ with neither $f_1$ nor $f_2$ constant.
Equivalently, this means that over the algebraic closure,
we can write $V=V_1\cup V_2$ nontrivially.
Since $V$ is nonsingular, this must be a disjoint union.
We get a contradiction, and hence a proof of Proposition \ref{agprop},
by showing the intersection number
$V_1\cdot V_2$ cannot be zero.
By the Riemann-Roch formula for the surface $Y$, we have
$\frac{1}{2}V\cdot V = {\rm Area}(\Gamma)$.
Letting $\Gamma_i$ denote the Newton polygon of $f_i$, for $i=1,2$,
we obtain $\frac{1}{2}V_1\cdot V_2 = \frac{1}{2}[ {\rm Area}(\Gamma) -
{\rm Area}(\Gamma_1) - {\rm Area}(\Gamma_2)] > 0$.
This positive quantity is the {\em mixed volume} of
$\Gamma_1$ and $\Gamma_2$; cf.\ \cite{f} or \cite{h}.
So $V$ is absolutely irreducible, and Proposition \ref{agprop} is proved.

\begin{proof}[Proof of upper bound in Theorem \ref{thm-toric}]
We use the notation of the proof of Lemma \ref{arithgenusformula}.
Fix $q$, and let $C$ be a curve of genus $g$ which embeds in a
toric surface $Y$ over $\Fq$.
Without loss of generality, we may assume $Y$ is nonsingular projective,
and has fixed point set (for the torus action) disjoint from $C$.
Since every rational point of $C$ is either in the torus $T^2$
or in one of the nontrivial intersections (\ref{modpu}),
the number of rational points on $C$ is at most $(q-1)^2+(q-1)v$,
where $v$ is the number of vertices of the Newton polygon of a defining
equation $f(x,y)$ for $C$ (in the coordinates of some toric chart).
So, it suffices to show that the minimum number $g(v)$ of
interior lattice points in a convex lattice $v$-gon satisfies
\begin{equation}
\label{boundgv}
g(v)\ge N\cdot v^3
\end{equation}
whenever $v\ge v_0$, for appropriate positive constants $N$ and $v_0$.

Arnol'd~\cite{A} showed that any convex lattice $v$-gon has area at
least $(1/8192)v^3$.  The desired bound (\ref{boundgv}) follows by
combining this, Pick's theorem, and the observation
that there is a $g$-minimal $v$-gon with no lattice points on the
boundary other than vertices (\cite{Ra,Si}: removing the triangle
bounded by two vertices and one interior edge point of a convex lattice
$v$-gon yields a $v$-gon with as many or fewer interior lattice points
and strictly smaller area).
\end{proof}


\section{Tame cyclic covers}
\label{sec-tame}

In this section we use cyclic, tamely ramified covers of $\Line$ to
produce curves with many points in every genus.  Let $\Nqtc$ denote the
maximal number of $\Fq$-rational points on any curve $C$ over $\Fq$ of
genus $g$ which admits a tame cyclic cover $C\to\Line$ over $\Fq$.
We will show

\begin{thm}
\label{thm-tame}
For any fixed $q$, we have
$$
\liminf_{g\to\infty} \frac{\Nqtc}{g (\log\log g)/(\log g)^3} >0.
$$
\end{thm}

At the end of
this section we discuss possible improvements of this result.
Our proof of Theorem~\ref{thm-tame} relies on the following lemma:

\begin{lemma}
\label{tamelemma}
Fix $q$ and $n>1$.  Let $\ell_1$ be the least prime not dividing $q$,
and let $\ell_2<\dots<\ell_n$ be primes which do not divide
$q(q-1)\ell_1$.  If $q$ is even, we assume further that
$\ell_i\equiv 7$ {\rm (mod}~$8)$ for some $i$.
Let $L=\prod_{i=1}^n \ell_i$.
Then, for any $g$ such that
\begin{equation}
\label{crt}
g>1-L+\frac{L}{2}\sum_{i=2}^n (\ell_i-1)^2,
\end{equation}
there is a curve $C/\Fq$ of genus $g$ and a tame cyclic
cover $C\to\Line$ (over $\Fq$) of degree $L$ in which
some $\Fq$-rational point of $\Line$ splits completely.
\end{lemma}

One can obtain $\liminf_{g\to\infty} \Nqtc (\log g)^4/(g\log\log g) > 0$
by taking $\ell_2$, $\dots$, $\ell_n$ to be the $n-1$ smallest
primes which do not divide $q(q-1)\ell_1$.
A slight modification to this choice of primes allows us to save a
factor of $\log g$ and thereby prove Theorem~\ref{thm-tame}.

\begin{proof}[Proof that Lemma~\ref{tamelemma} implies Theorem~\ref{thm-tame}]
Fix $q$, and let $\ell_1$ be the least prime not dividing $q$.
For any $g\geq 0$, define $x_g$ to be the least integer such that
\begin{equation}
\label{xg}
 g-1 \leq \frac{\ell_1}{2}\cdot\Biggl(\prod_{\substack{\ell\leq x_g \\
 \ell\nmid q(q-1)\ell_1 \\ \ell\text{ prime}}}
 \ell\Biggr)\cdot\Biggl(-2+\sum_{\substack{\ell\leq x_g \\
 \ell\nmid q(q-1)\ell_1 \\ \ell\text{ prime}}} (\ell-1)^2\Biggr).
\end{equation}
Note that the right side might be only slightly larger than the left,
or it might be larger by a factor of as much as (slightly more than) $x_g$.
We will modify
the set of primes under consideration in order to find an analogous
product which is only slightly smaller than $g-1$; we do this by
replacing two (suitably chosen) primes $p_1$ and $p_2$ by a third prime
$p_3$.

We now define $p_1$, $p_2$, and $p_3$; each definition makes sense
for $g$ sufficiently large.  Let $p_1$ be the smallest prime whose
removal from the right side of (\ref{xg}) would reverse the inequality;
that is, $p_1$ is the least prime such that $p_1\nmid q(q-1)\ell_1$ and
\begin{equation*}
 g-1 > \frac{\ell_1}{2p_1}\cdot\Bigl(\prod_\ell
 \ell\Bigr)\cdot\Bigl(-2-(p_1-1)^2+\sum_\ell (\ell-1)^2\Bigr).
\end{equation*}
(Here, and in the remainder of this proof, any sum or product indexed by
 $\ell$ is understood to be taken over all primes $\ell$ such that
 $\ell\leq x_g$ and $\ell\nmid q(q-1)\ell_1$.)
Let $p_2$ be the largest prime such that $p_2\leq x_g$ and
$p_2\nmid q(q-1)\ell_1p_1$.  Let $p_3$ be the largest prime
such that
\begin{equation}
\label{noc}
\begin{split}
g-1>&\frac{\ell_1p_3}{2p_1p_2}\cdot\Bigl(\prod_\ell
 \ell\Bigr)\cdot \\&\Bigl(-2-(p_1-1)^2-(p_2-1)^2+(p_3-1)^2+
 \sum_\ell (\ell-1)^2\Bigr).
\end{split}
\end{equation}

We apply Lemma~\ref{tamelemma} to the set of primes
\begin{equation*}
\{\ell_1,\dots,\ell_n\} := \{\ell\leq x_g\colon \ell\nmid q(q-1)p_1p_2,
  \text{ $\ell$ prime}\}  \cup \{p_3\} \cup \{\ell_1\}.
\end{equation*}
(The hypotheses of Lemma~\ref{tamelemma} are satisfied when $g$ is sufficiently
large, since then $p_3>\max\{q,\ell_1\}$ and also there will be
some $i$ for which $\ell_i\equiv 7$ (mod~8).)
It follows that $\Nqtc\geq \prod_{i=1}^n\ell_i=(\prod_\ell\ell)\cdot
\ell_1 p_3/(p_1 p_2)$.
It remains only to determine the asymptotics when $g\to\infty$.

We can rewrite (\ref{xg}) as
\begin{equation}
\label{xg2}
\log (g-1)\leq\log(\ell_1/2)+
 \Bigl(\sum_\ell \log\ell\Bigr) + \log\Bigl(-2+\sum_\ell (\ell-1)^2\Bigr).
\end{equation}
The Prime Number Theorem implies that the right hand side of (\ref{xg2})
is asymptotic to $x_g$ as $x_g\to\infty$, so $x_g\sim\log g$.
Note that $\sum_{\ell} (\ell-1)^2$ is asymptotic to $x_g^3/(3\log x_g)$;
this is much larger than $p_1$, $p_2$, and $p_3$, since
$p_1,p_2\leq x_g$ and (for $x_g$ large) $p_3<5qx_g$.  Since the left and
right sides of (\ref{noc}) are asymptotic to each other as $x_g\to\infty$,
it follows that
$$ \Nqtc\geq\frac{\ell_1 p_3}{p_1p_2} \cdot \Bigl(\prod_\ell \ell\Bigr) \sim
  \frac{2g}{x_g^3/(3\log x_g)} \sim \frac{6g(\log\log g)}{(\log g)^3}.
$$
This completes the proof.
\end{proof}

Our proof of Lemma~\ref{tamelemma} uses the following existence result:

\begin{lemma}
\label{tamecftlemma}
Let $P_1, \dots, P_n,\infty$ be distinct places on $\Line/\Fq$, with
degrees $d_1, \dots, d_n,1$.  Let $\ell_1, \dots, \ell_n$ be positive
integers such that $q^{d_i}\equiv 1$ {\rm (mod}~$\ell_i(q-1))$ for each $i$.
There exists a tame abelian cover $\phi\colon C\to\Line$ (over $\Fq$)
of degree $L :=\prod_{i=1}^n \ell_i$ such that $\infty$ splits
completely in $C$ and the genus $g$ of $C$ satisfies
\begin{equation}
\label{tamecfteq}
2g-2 = -2L + L\sum_{i=1}^n \frac{\ell_i-1}{\ell_i}d_i.
\end{equation}
Moreover, if the $\ell_i$ are pairwise coprime then $\phi$
can be chosen to be cyclic.
\end{lemma}

\begin{proof}[Proof of Lemma~\ref{tamecftlemma}]
Let $\widehat{\phi}\colon \widehat{C}\to\Line$ be the maximal
tamely ramified abelian cover of $\Line/\Fq$ which is unramified
outside $\{P_1,\dots,P_n\}$ and in which $\infty$ splits completely.
By Class Field Theory for $\Line$, the Galois group $\widehat{G}$
of $\widehat{\phi}$ fits in a short exact sequence
$$
\begin{CD}
        1 @>>> \F_q^\times @>\Delta>> \prod_{i=1}^{n}
             \F_{q^{d_i}}^\times @>>> \widehat{G}  @>>> 1,
\end{CD}
$$
where $\Delta$ is the diagonal embedding into the product.
Moreover, the inertia group over the place $P_i$ is the
image of $\F_{q^{d_i}}^\times$ in $\widehat{G}$.

Since $\ell_i$ divides $(q^{d_i}-1)/(q-1)$, the group
$G:= \prod_{i=1}^n \Z/\ell_i\Z$ is a quotient of $\widehat{G}$;
let $\phi\colon C\to\Line$ be the corresponding cover.
Then $\phi$ is a tame abelian cover with Galois group $G$
in which $\infty$ splits completely and in which all places
outside $\{P_1,\dots,P_n\}$ are unramified.  Moreover, the inertia group
over the place $P_i$ is $\Z/\ell_i\Z$.  Hence the genus $g$ of $C$
satisfies (\ref{tamecfteq}).

Finally, if the $\ell_i$ are pairwise coprime then $G\cong\Z/L\Z$
is cyclic.
\end{proof}

\begin{proof}[Proof of Lemma~\ref{tamelemma}]
Assume $q$ is odd, so $\ell_1=2$.
Let $s_1,\dots,s_n$ be positive integers which
will be specified later.
Put $r_1=2$ and $r_i=\ell_i-1$
for $i>1$; set $d_i = r_i s_i$ for all $i$.
Note that $q^{d_i}\equiv 1$ (mod~$\ell_i(q-1)$) for all $i$.

An easy count shows that, for
any $d>0$, there are at least $d$ finite places on $\Line$ of degree $d$;
since $d_i\geq i$, it follows that we can choose distinct finite places
$P_1,\dots,P_n$ on $\Line$ with degrees $d_1,\dots,d_n$.
%
%
Lemma~\ref{tamecftlemma} yields a tame cyclic cover $C\to\Line$ (over $\Fq$)
of degree $L :=\prod_{i=1}^n \ell_i$ such that some $\Fq$-rational point of
$\Line$ splits completely in $C$ and the genus $\gt$ of $C$ satisfies
\begin{equation*}
2\gt-2 = -2L + L\sum_{i=1}^n \frac{\ell_i-1}{\ell_i}d_i.
\end{equation*}

We must show that, for any $g$ satisfying (\ref{crt}), we can choose
the $s_i$ so that $\gt=g$.  Pick any $g$ satisfying (\ref{crt}).
Rewrite the expression for
$\gt$ as
\begin{equation}
\label{tildeeq}
\gt-1+L = \frac{L}2 s_1 + \sum_{i=2}^n
\frac{L}{2\ell_i} (\ell_i-1)^2 s_i.
\end{equation}
Note that, for $i>1$, we have $\gt-1\equiv s_i L/(2\ell_i)$
(mod~$\ell_i$).  For each $i>1$, let $s_i$ be the unique integer
such that $1\leq s_i\leq \ell_i$ and
$g-1\equiv s_i L/(2\ell_i)$ (mod~$\ell_i$).  Then
$g-1\equiv\gt-1$ (mod~$L/2$), so (\ref{tildeeq}) implies there is
a unique integer $s_1$ for which $\gt=g$.  It remains to show $s_1>0$;
this follows from (\ref{crt}), since
$$
g-1+L > \frac L{2}\sum_{i=2}^n (\ell_i-1)^2
       \geq \sum_{i=2}^n \frac{L}{2\ell_i}(\ell_i-1)^2 s_i
         = \gt-1+L-\frac{L}2 s_1.
$$

Finally, we indicate how the argument must be modified to
handle the case of even $q$.  It suffices to prove the result
for $q=2$.  In this case, let $i_0$ satisfy
$1<i_0\leq n$ and $\ell_{i_0}\equiv 7$ (mod~8).
Define $r_i$ as above for $i\neq i_0$, and let $r_{i_0}=
(\ell_{i_0}-1)/2$.  As above,
%
%
%
%
there is a tame cyclic cover $C\to\Line$ (over $\F_2$) of
degree $L:=\prod_{i=1}^n \ell_i$ such that some $\F_2$-rational point splits
completely in $C$ and the genus $\gt$ of $C$
satisfies~(\ref{tildeeq}).  It remains to choose the $s_i$
so that $\gt=g$.  For $i\notin\{1,i_0\}$, choose $s_i$
as above; for $i=i_0$, let $s_i$ satisfy $1\leq s_i\leq 2\ell_i$
and $g-1+L\equiv s_i (L/\ell_i)((\ell_i-1)/2)^2$ (mod~$2\ell_i$).
Then (\ref{crt}) implies there is a unique positive integer $s_1$
such that $\gt=g$, which completes the proof of Lemma~\ref{tamelemma},
and thus the proof of Theorem~\ref{thm-tame}.
%
%
\end{proof}


\begin{remark}
Theorem~\ref{thm-tame} implies that $\Nqtc > C_q\cdot g(\log\log g)/(\log g)^3$
for $g$ sufficiently large, where $C_q$ is a positive constant depending only on $q$.  We
do not know whether this result can be improved.  In the opposite direction,
we now determine the qualitatively best possible {\em upper} bound on $\Nqtc$.
It was shown in~\cite{FPS} that $\Nqtc < D_q\cdot g/\log g$ for $g>1$, where $D_q$
is a positive constant depending only on $q$.  The following
result shows that this upper bound is best possible, by exhibiting
infinitely many $g$ (for each $q$) such that $\Nqtc>g/\log g$.
\end{remark}

\begin{prop}
For any fixed $q$, there are infinitely many $g$ for which we have
$\Nqtc>(2\log{q}) g/\log{g}$.
\end{prop}

\begin{proof}
Fix $q$.  For any $e\geq 4$, put $d=(q^e-1)/(q-1)$ and let $P,\infty$
be places of $\Line$ of degrees $e$ and $1$.  Let $C\to\Line$ be the
maximal tame abelian cover in which all places besides $P$ are unramified
and in which $\infty$ splits completely.  Then $C\to\Line$ is a cyclic
cover of degree $d$ in which $P$ is totally ramified.  Moreover,
since $\infty$ splits completely, the cover $C\to\Line$ is defined
over $\Fq$ and $\#C(\Fq)\geq d$.  The genus of $C$ is $g=(d-1)(e/2-1)$.
Finally,
\begin{equation*}
\frac{g}{d} < \frac{e-1}{2} < \frac{\log(d-1)}{2\log q} \le \frac{\log g}{2\log q},
\end{equation*}
so $d>(2\log{q})g/\log{g}$ and thus the proof is complete.
\end{proof}





\begin{thebibliography}{99}

\newcommand{\au}[1]{{#1},}
\newcommand{\ti}[1]{{#1},}
\newcommand{\jo}[1]{\textit{#1}}
\newcommand{\vo}[1]{\textbf{#1}}
\newcommand{\yr}[1]{(#1),}
\newcommand{\pp}[1]{#1.}
\newcommand{\bk}[1]{\textit{#1},}
\newcommand{\inbk}[1]{in: \bk{#1}}
\newcommand{\inbkpp}[1]{pp.~#1.}
\newcommand{\xxx}[1]{{arXiv:#1}}

\bibitem{A}
\au{Arnold, V. I.}
\ti{Statistics of integral convex polygons}
\jo{Funktsional. Anal. i Prilozhen.}
\vo{14}
\yr{1980}
\pp{1--3}
[\jo{Funct. Anal. Appl.}
\vo{14}
\yr{1980}
\pp{79--81}]



\bibitem{CWZ}
\au{Csirik, J., Wetherell, J. and Zieve, M.}
\ti{On the genera of $X_0(N)$}
preprint, 2000,
\xxx{math.NT/0006096}.

\bibitem{DV}
\au{Drinfeld, V. G. and Vladut, S. G.}
\ti{The number of points of an algebraic curve}
\jo{Funktsional. Anal. i Prilozhen.}
\vo{17}
\yr{1983}
\pp{68--69}
[\jo{Funct. Anal. Appl.}
\vo{17}
\yr{1983}
\pp{53--54}]

\bibitem{El98}
\au{Elkies, N.}
\ti{Explicit modular towers}
\inbk{Proceedings of the Thirty-Fifth Annual Allerton Conference on
    Communication, Control and Computing}
(ed. T. Basar and A. Vardy),
Univ. of Illinois at Urbana-Champaign, 1998,
\inbkpp{23--32}

\bibitem{El3ecm}
\au{Elkies, N.}
\ti{Explicit towers of Drinfeld modular curves}
\inbk{Proceedings of the Third European Congress of Mathematics
(Barcelona, 2000)}
to appear.
\xxx{math.NT/0005140}

\bibitem{phewkz}
\au{Elkies, N., Howe, E., Kresch, A., Poonen, B., Wetherell, J. and Zieve, M.}
\ti{Curves of every genus with many points, II: On a question of Serre}
preprint, 2000.

\bibitem{FKV}
\au{Frey, G., Kani, E. and V\"olklein, H.}
\ti{Curves with infinite $K$-rational geometric fundamental group}
\inbk{Aspects of Galois Theory (Gainesville, Fla., 1996)}  
London Math. Soc. Lect. Note Ser. {\bf 256},
Cambridge Univ. Press, Cambridge, 1999,
\inbkpp{85--118}

\bibitem{FPS}
\au{Frey, G., Perret, M. and Stichtenoth, H.}
\ti{On the different of abelian extensions of global fields}
\inbk{Coding Theory and Algebraic Geometry (Luminy, 1991)}  
Lect. Notes in Math. {\bf 1518},
Springer-Verlag, New York, 1992,
\inbkpp{26--32}

\bibitem{f}
\au{Fulton, W.}
\bk{Introduction to Toric Varieties}
Ann. of Math. Stud., {\bf 131},
Princeton Univ. Press, Princeton, 1993.

\bibitem{GS0}
\au{Garcia, A. and Stichtenoth, H.}
\ti{Elementary abelian $p$-extensions of algebraic function fields}
\jo{Manuscripta Math.}
\vo{72}
\yr{1991}
\pp{67--79}

\bibitem{GS1}
\au{Garcia, A. and Stichtenoth, H.}
\ti{A tower of Artin-Schreier extensions of function fields attaining
    the Drinfeld-Vladut bound}
\jo{Invent. Math.}
\vo{121}
\yr{1995}
\pp{211--222}

\bibitem{GS3}
\au{Garcia, A. and Stichtenoth, H.}
\ti{On the asymptotic behavior of some towers of function fields over
    finite fields}
\jo{J. Number Theory}
\vo{61}
\yr{1996}
\pp{248--273}

\bibitem{GST}
\au{Garcia, A., Stichtenoth, H. and Thomas, M.}
\ti{On towers and composita of towers of function fields over finite fields}
\jo{Finite Fields Appl.}
\vo{3}
\yr{1997}
\pp{257--274}

\bibitem{GV}
\au{van der Geer, G. and van der Vlugt, M.}
\ti{Tables of curves with many points}
\jo{Math. Comp.}
\vo{69}
\yr{2000}
\pp{797--810}
Updates on \verb+http://www.science.uva.nl/~geer/+

\bibitem{GVJCT}
\au{van der Geer, G. and van der Vlugt, M.}
\ti{Fiber products of Artin-Schreier curves and generalized Hamming weights
of codes}
\jo{J. Combin. Theory Ser. A}
\vo{70}
\yr{1995}
\pp{337--348}


\bibitem{h}
\au{Hovanskii, A. G.}
\ti{Newton polyhedra, and the genus of complete intersections}
\jo{Funktsional. Anal. i Prilozhen.}
\vo{12}
\yr{1978}
\pp{51--61}
[\jo{Funct. Anal. Appl.}
\vo{12}
\yr{1978}
\pp{38--46}]

\bibitem{Ih-66}
\au{Y. Ihara}
\ti{Algebraic curves mod $\mathfrak{p}$ and arithmetic groups}
\inbk{Algebraic Groups and Discontinuous Subgroups (Boulder, Colo., 1965)}
Proc. Sympos. Pure Math. {\bf 9},
Amer. Math. Soc., Providence, 1966,
\inbkpp{265--271}

\bibitem{Ih-69}
\au{Y. Ihara}
\bk{On Congruence Monodromy Problems, Vol. 2}
Department of Mathematics, Univ. of Tokyo, 1969.

\bibitem{Ih-75}
\au{Y. Ihara}
\ti{On modular curves over finite fields}
\inbk{Discrete Subgroups of Lie Groups and Applications to Moduli (Internat.
Colloq., Bombay, 1973)}
Oxford Univ. Press, Bombay, 1975,
\inbkpp{161--202}

\bibitem{Ih-79}
\au{Y. Ihara}
\ti{Congruence relations and Shim\=ura curves}
\inbk{Automorphic Forms, Representations, and $L$-functions
(Corvallis, Ore., 1977)}
Proc. Sympos. Pure Math. {\bf 33}, Part 2,
Amer. Math. Soc., Providence, 1979,
\inbkpp{291--311}

\bibitem{Ih}
\au{Ihara, Y.}
\ti{Some remarks on the number of rational points of algebraic curves
    over finite fields}
\jo{J. Fac. Sci. Univ. Tokyo}
\vo{28}
\yr{1981}
\pp{721--724}

\bibitem{MV}
\au{Manin, Y. I. and Vladut, S. G.}
\ti{Linear codes and modular curves}
\jo{Itogi Nauki i Tekhniki}
\vo{25}
\yr{1984}
\pp{209--257}
[\jo{J. Soviet Math.}
\vo{30}
\yr{1985}
\pp{2611--2643}]

\bibitem{NX0}
\au{Niederreiter, H. and Xing, C. P.}
\ti{Algebraic curves over finite fields with many rational points}
\inbk{Number Theory (Eger, 1996)}
W. de Gruyter, Berlin,
\yr{1998}
\inbkpp{423--443}

\bibitem{NX}
\au{Niederreiter, H. and Xing, C. P.}
\ti{Towers of global function fields with asymptotically many rational
    places and an improvement on the Gilbert-Varshamov bound}
\jo{Math. Nachr.}
\vo{195}
\yr{1998}
\pp{171--186}

\bibitem{Ra}
\au{Rabinowitz, S.}
\ti{On the number of lattice points inside a convex lattice $n$-gon}
\jo{Congr. Numer.}
\vo{73}
\yr{1990}
\pp{99--124}

\bibitem{Se}
\au{Serre, J.-P.}
\ti{Sur le nombre des points rationnels d'une courbe alg\'ebrique sur un
    corps fini}
\jo{C. R. Acad. Sci. Paris}
\vo{296}
\yr{1983}
\pp{397--402; = \OE uvres [128]}

\bibitem{Se2}
\au{Serre, J.-P.}
\ti{Nombres de points des courbes alg\'ebriques sur $\Fq$}
\jo{S\'em. Th\'eor. de Nombres Bordeaux}
\yr{1982--1983}
exp. 22; = \OE uvres [129].

\bibitem{Se3}
\au{Serre, J.-P.}
\ti{R\'esum\'e des cours de 1983--1984}
\inbk{Ann. Coll\`ege de France} 1984,
pp.~79--83; = \OE uvres [132].

\bibitem{Se4}
\au{Serre, J.-P.}
\ti{Rational points on curves over finite fields}
unpublished lecture notes by F. Gouv\^ea, Harvard Univ., 1985.


\bibitem{Si}
\au{Simpson, R.}
\ti{Convex lattice polygons of minimum area}
\jo{Bull. Austral. Math. Soc.}
\vo{42}
\yr{1990}
\pp{353--367}

\bibitem{TVZ}
\au{Tsfasman, M. A., Vladut, S. G. and Zink, Th.}
\ti{Modular curves, Shimura curves, and Goppa codes, better than
Varshamov-Gilbert bound}
\jo{Math. Nachr.}
\vo{109}
\yr{1982}
\pp{21--28}

\bibitem{Zi}
\au{Zink, Th.}
\ti{Degeneration of Shimura surfaces and a problem in coding theory}
\inbk{Fundamentals of Computation Theory (Cottbus, 1985)}
Lect. Notes in Comput. Sci. {\bf 199},
Springer-Verlag, New York, 1985,
\inbkpp{503--511}

\end{thebibliography}
\end{document}